\documentclass[draft,12pt,leqno]{amsart}

 \usepackage{amsmath,amssymb,amscd,amsthm}
 \usepackage{latexsym}

 \newtheorem{theorem}{Theorem}
 \newtheorem{lemma}{Lemma}
 \newtheorem{proposition}{Proposition}
 
 \newtheorem{question}{Question}

 \newcommand{\q}{\quad}
 \newcommand{\qq}{\quad\quad}

 \newcommand{\norm}[2]{{\left\| #1 \right\|}_{#2}}

 \newcommand{\dD}{\partial D}
 \newcommand{\ds}{\displaystyle}

 \newcommand{\al}{\alpha}
 
 \newcommand{\ga}{\gamma}
 \newcommand{\Ga}{\Gamma}
 \newcommand{\de}{\delta}
 \newcommand{\De}{\Delta}
 \newcommand{\ve}{\varepsilon}

 \newcommand{\la}{\lambda}

 \newcommand{\si}{\sigma}
 \newcommand{\Si}{\Sigma}
 \newcommand{\vp}{\varphi}
 
 \newcommand{\Om}{\Omega}

 \newcommand{\rn}{{\mathbb R}^n}
 \newcommand{\rone}{\mathbb R^1}
 \newcommand{\rtwo}{\mathbb R^2}

 \newcommand{\lp}{L^{p}}

 \newcommand{\hpo}{H^p(\dD)}
 
 \newcommand{\ca}{C^\al(\dD)}

 \newcommand{\cs}{\mathcal S}
 \newcommand{\cc}{\mathcal C}

 \newcommand{\ck}{\mathcal K}

 \newcommand{\pv}{\textup{p.v.}\,}
 
 \newcommand{\intl}{\int\limits}

 \newcommand{\suml}{\sum\limits}
 \newcommand{\supl}{\sup\limits}
 \newcommand{\f}{\displaystyle\frac}
 \newcommand{\p}{\partial}
 \newcommand{\pp}[2]{\f{\p #1}{\p #2}}
 
 \newcommand{\R}{\operatorname{R}}
 \newcommand{\dpr}[2]{\langle #1, #2\rangle}
 \newcommand{\di}{\textup{div}}
 \newcommand{\arr}[4]{\left(\begin{array}{cc} #1 & #2\\#3 & #4\end{array}\right)}

\begin{document}

 \title[Optimal solvability for the Dirichlet and  Neumann problem]
 {Optimal solvability for the Dirichlet and  Neumann problems in
 dimension two}

 \author[Stefanov, Verchota]
{A. Stefanov, G.C. Verchota
\address{Atanas Stefanov, LGRT 1342, Department of Mathematics and Statitics,
University of Massachusetts, Amherst, MA 01003, USA} \\
\email{stefanov@math.umass.edu}
\address{Gregory Verchota, 215, Carnegie Hall,
Department of Mathematics, Syracuse University, Syracuse, NY
13088, USA}}

\begin{abstract}
We show  existence and uniqueness for the solutions of the
regularity and the Neumann problems for harmonic functions on
Lipschitz domains with data in the Hardy spaces $H^p_1(\dD)$($H^{p}(\partial
D)$), $\displaystyle p>\f{2}{3}-\ve$, where $D\subset \rtwo$ and
$\ve$ is a (small) number depending on the Lipschitz nature of
$D$.
This in turn implies that solutions to the Dirichlet
problem with data in the H\"older class $C^{1/2+\ve}(\partial D)$
are themselves in $C^{1/2+\ve}(\bar{D})$.  Both of these results
are sharp.
In fact, we prove a more general statement
regarding the $H^p$ solvability for divergence form elliptic
equations with bounded measurable coefficients.
 
We also provide  $H^{2/3-\ve}$ and $C^{1/2+\ve}$
solvability result for the regularity and Dirichlet problem 
for the biharmonic equation on Lipschitz domains.
\end{abstract}

\maketitle
\date{}

 \section{Introduction and Main Results}

 In this paper we study the Dirichlet and Neumann problems for
 harmonic functions on Lipschitz domains and their biharmonic counterparts. 
 More precisely let $X, Y,
 Z$ be function spaces on the boundary $\dD$ of $D$. Then $$
 (D_X)\left| \begin{array}{cc} \De u &=0\\ u|_{\dD} &=f \\ M(u)
 &\in X
 \end{array}\right.$$ is the Dirichlet problem with underlying
 space $X$, and $$(N_Y)\left| \begin{array}{cc} \De u &=0\\ \left.
 \ds\pp{u}{N} \right|_{\dD} &=g
 \\ M(\nabla u) &\in Y
 \end{array}\right.$$  $$ (R_Z)\left| \begin{array}{cc}
 \De u &=0\\ u|_{\dD} &=h \\ M(\nabla u) &\in Z
 \end{array}\right.$$
 are the Neumann and  regularity problems. Here $N$ is the outer
 normal vector to $D$, $\ds\pp{u}{N}=\dpr{N}{\nabla u}$ and $M(u)$
 is the usual non-tangential maximal function of $u$.

  In this setting the canonical choices are $X=Y=\lp(\dD)$,
$Z=\lp_1(\dD)$- the
  space of functions with one distributional derivative in
 $\lp(\dD)$. In the sequel, we will slightly abuse notations by
 using $D_p$ instead of $D_{\lp}$, $N_p$ instead of $N_{\lp}$
 etc.
\subsection{Harmonic functions}
  We state now the
 classical results related to the $\lp$ theory.\\ \\ {\bf Theorem}
 \cite{Dahlberg1, Verchota, Kenig0, Kenig} Let $D\subset \rn$ be a connected Lipschitz
 domain. Then
 \begin{enumerate}
 \item there exists an $\ve=\ve(D)>0$ such that for $2-\ve<p\leq
 \infty$ and $f\in\lp(\dD)$ there is an unique solution to $D_p$.
 Moreover, there is the {\it apriori} estimate \\
 $\norm{M(u)}{\lp(\dD)}\leq C \norm{f}{\lp(\dD)}$.
 \item there exists an $\ve=\ve(D)>0$ such that for $1<p<2+\ve$ and
 $g\in\lp(\dD), \int_{\dD} g d\si =0$ there is an unique (up to a
 constant) solution to $N_p$. There is the {\it apriori} estimate
 \\ $\norm{M(\nabla u)}{\lp(\dD)}\leq C \norm{g}{\lp(\dD)}$.
 \item there exists an $\ve=\ve(D)>0$ such that for $1<p<2+\ve$ and
 $h\in\lp_1(\dD)$ there is an unique solution to $R_p$. There is
 the {\it apriori} estimate
 \\ $\norm{M(\nabla u)}{\lp(\dD)}\leq C \norm{h}{\lp_1(\dD)}$.
 \end{enumerate}
This theorem summarizes the results in \cite{Dahlberg1}, but 
some earlier version and ideas originated in \cite{Verchota}. Actually,the 
$L^p$ theory described above is a consequence of the duality between the 
Dirichlet and regularuty problems, the $L^2$ solvability for all three 
problems and the following endpoint result due to Dahlberg and Kenig.
\begin{theorem}[Dahlberg-Kenig]
 \label{theo:th1} Let $D\subset \rn$ be a connected star-like
 Lipschitz domain. Then\begin{enumerate}
 \item $$ (N_1)\left| \begin{array}{cc} \De u
 &=0\\ \left. \ds\pp{u}{N} \right|_{\dD} &=f
 \\ M(\nabla u) &\in L^1
 \end{array}\right.$$
 is uniquely solvable provided $f\in H^1(\dD)$ and $\norm{M(\nabla
 u)}{L^1(\dD)}\leq C \norm{f}{H^1(\dD)}$.
 \item Given $f\in H^1_1(\dD)$ there exists an unique solution to
 $$ (R_1)\left| \begin{array}{cc} \De u
 &=0\\  u|_{\dD} &=f
 \\ M(\nabla u) &\in L^1
 \end{array}\right.$$ Moreover
 $\norm{M(\nabla u)}{L^1(\dD)}\leq C \norm{f}{H^1_1(\dD)}$.
 \end{enumerate}
 \end{theorem}
 As a corollary to this result, one  proves (weak) maximum
 principle, solvability for $BMO$ data etc. We refer to
 \cite{Dahlberg1} for excellent treatise of these questions.

 Recently, Brown \cite{Brown} was able to extend Theorem
 \ref{theo:th1} to show that there exists $\ve=\ve(D)$ , such that
 for $1-\ve<p<1$ the Neumann problem $N_p$  is still uniquely
 solvable with the usual estimates $\ds\norm{M(\nabla u)}{p}\leq
 C\norm{f}{H^p}$. This result has the interesting corollary that
 the double layer potential is invertible operator on the H\"older
 space $C^\al(\dD)$ for $\al$ close to zero and thus we have a
 representation formula for the solutions of the Dirichlet problem
 with $C^\al$ data. This raises the following natural question, see
 Question 3.2.10 in  \cite{Kenig}.
 \begin{question}
 \label{question1}
  Are $R_p$ and $N_p$ solvable for $p$
 significantly below one? What does that imply for  solutions of
 the Dirichlet problem with $C^\al$ data for $\al$ significantly
 above zero?
 \end{question}
 The purpose of this paper is to establish the optimal $p$ range
 for solvability of both $R_p$ and $N_p$ in dimension two. That is
 our Theorem \ref{theo:th2} below.  Let us remark only that known
 counterexamples in dimensions bigger than two imply that  the
 Neumann problem $N_p$ may not be solvable for $p<1-\ve(D)$, i.e. for
 fixed $p<1$, there exists a Lipschitz domain $D$ such that
 $N_{p}(D)$ is not uniquely solvable.
 \begin{theorem}
 \label{theo:th2} Let $D\subset \rtwo$ be a star-like Lipschitz
 domain with connected boundary. There exists $\ve=\ve(D)$, such
 that for $2/3-\ve<p<1$ and $0<\al<1/2+\ve$
 \begin{enumerate}
 \item The Neumann problem $$(N_p)\left| \begin{array}{cc}& \De u
 =0\\ & \ds\pp{u}{N} \left|_{\dD} =f\in H^p(\dD)\right.
 \end{array}\right.$$ is uniquely solvable and $\norm{M(\nabla
u)}{L^p(\dD)}\leq
 C\norm{f}{H^p(\dD)}$.
 \item The regularity problem
 $$ (R_p)\left| \begin{array}{cc} &\De u =0\\ & u|_{\dD} =f\in
 H^p_1(\dD)
 \end{array}\right.$$ has unique solution and
 $\norm{M(\nabla u)}{L^p(\dD)}\leq C \norm{f}{H^p_1(\dD)}$.
 \item The Dirichlet problem
 $$ (D_\al)\left| \begin{array}{cc} &\De u =0\\ & u|_{\dD} =f\in
 C^\al(\dD)
 \end{array}\right.$$ has unique solution and
 $\norm{u}{C^\al(\bar{D})}\leq C \norm{f}{C^\al(\dD)}$.
 \end{enumerate}
 Moreover, the ranges $2/3-\ve<p$ and $\al<1/2+\ve$ are sharp.
 \end{theorem}
 In fact, we consider more general divergence form elliptic
 equations in the form \\ $\di(A(\nabla u)) = 0$, where $A$ is a
 symmetric, elliptic matrix with real-valued bounded measurable
coefficients.
  We prove the following theorem.
 \begin{theorem}
 \label{theo:main} Let $\ds A(x,t)=A(x)=\left(\begin{array}{cc}
 a(x) & b(x)
 \\b(x) & c(x)\end{array}\right) $ be a real, symmetric, uniformly
 elliptic
 matrix with bounded and measurable coefficients,  independent of
 the time variable. Then there exists  $\ve=\ve(D)>0$ such that for
 $2/3-\ve<p<2+\ve$
 \begin{enumerate}
 \item the Neumann problem in the upper half-space
 $\rtwo_+=\{(x,t): t>0\}$ $$ (N_p)\left| \begin{array}{cc} &\di(A
 \nabla u) =0\ \  \textup{for}\ t>0,
 \\ & A\nabla u(x,0) \cdot (0,-1)| =f\in H^p(\rone)
 \end{array}\right.$$
 has unique solution and $\norm{M(\nabla u)}{\lp}\leq
 C\norm{f}{H^p}$.
 \item The regularity problem in the upper half-space
 $\rtwo_+$ $$ (R_p)\left| \begin{array}{cc} &\di(A \nabla u) =0 \ \
 \textup{for}\ t>0,
 \\ & u(x,0) =h\in H^p_1(\rone)
 \end{array}\right.$$
 has unique solution and $\norm{M(\nabla u)}{\lp}\leq
 C\norm{h}{H^p_1}$.
 \end{enumerate}
 Moreover, the range $2/3-\ve<p$ is sharp.
 \end{theorem}

 {\bf Remarks}
 \begin{itemize}
\item For divergence form equations
 $\di (A\nabla u)=0$, $A=A(x,t)$ one cannot expect solvability even
 for $D_2$ or $N_2$. Indeed, counterexamples show that unless we
 require {\it radial independence} for such a problem in the unit
 ball, we may encounter non-uniqueness for $N_2$, see \cite{Kenig1}
 and \cite{Kenig}, p. 63. \\ We will however work in the upper
 half-space instead of the unit ball. These two problems are not so
 much different. In fact, the appropriate assumption in the
 upper-half space is {\it time independence} (see \cite{Kenig},
 p.68 for a relevant discussion) and so our theorem \ref{theo:main}
 is formulated in that fashion. Let us remark only, that the
 problems $D_2$, $R_2$ and $N_2$ are all solvable for matrices
 $A=A(x)$ with time independent coefficients (\cite{Kenig1}). We make
use
 of these facts later on in our  proofs.
 \item The restriction to the upper-half space  in theorem
\ref{theo:main}
 is just for technical reasons. In fact, one can state the theorem
 for a general Lipschitz domain in $\rtwo$. The following argument
 shows that for the Dirichlet (regularity) problem.
\end{itemize}

 Let $D$ be the domain above the Lipschitz graph $t=\varphi(x)$
 and $u$ solves the Dirichlet problem $\ds\di(A(x)\nabla u)=0$,
 $\ds u(x,\varphi(x))=f(x)$. Define $\Phi(x,t)=(x,t-\varphi(x))$
 and set $\widetilde{u}(\Phi(x,t))=u(x,t)$. It is not
 difficult to check that $\widetilde{u}:\rtwo_+\to \rone$ is a
 solution to  $\ds\di(\widetilde{A}\nabla \widetilde{u})=0$, $\ds
 \widetilde{u}(x,0)=f $, where $$
 \widetilde{A}(x)=\arr{1}{0}{-\varphi'(x)}{1}A(x)
 \arr{1}{-\varphi'(x)}{0}{1}.$$  In particular, we have shown that
 Theorem \ref{theo:main} implies parts  one and two of
 Theorem \ref{theo:th2}.

 Unfortunately, at this moment we cannot claim part three of our
 Theorem \ref{theo:th2} for general divergence form elliptic
 equations with time independent coefficients. Our proof for
 harmonic functions is based on Brown's duality
 technique for the double-layer potential, which does not seem to
 generalize in the setting of Theorem \ref{theo:main}. Thus we pose
 the following:
 \begin{question}
 Assume that $A=A(x)$ is a real, symmetric elliptic matrix. Prove
 that the Dirichlet problem in the upper half-space $\rtwo_{+}$ $$
 (D_\al)\left|
 \begin{array}{cc} &\di(A \nabla u) =0 \ \ \textup{for} \ t>0,
 \\ & u |_{\dD} =f\in C^\al(\rone)
 \end{array}\right.$$
 is solvable with  $\norm{u}{C^\al(\rtwo_+)}\leq
 C\norm{f}{C^\al(\rone)}$ for $\al<1/2+\ve(A)$.
\end{question}
We now state the following result, which gives a connection between 
the Neumann Hardy spaces and the usual atomic Hardy spaces. This is an 
extension of  Theorem 2.3.18 in \cite{Kenig}, for the case $p<1$.
 \begin{theorem}
\label{prop:five}
Let $D\subset \rtwo$ is a domain above Lipschitz graph and $u$ is 
a harmonic function on $D$ that satisfies $u|_{\dD}\in H^p_1(\dD)$, 
$2/3-\ve<p\leq 1$. Then  
$$\norm{\f{\p u}{\p N}}{H^p(\dD)}\lesssim \norm{u}{H^1_p(\dD)}.$$
Conversely, given $f\in H^p(\dD)$, there exists a harmonic function $u$,
such that $\f{\p u}{\p N}=f$ and 
$$
\norm{u}{H^p_1(\dD)}\lesssim \norm{f}{H^p(\dD)}.
$$
\end{theorem}

\subsection{Biharmonic functions}
For the biharmonic equation, we consider the Dirichlet problem 
 $$BD_p\left|\begin{array}{cc}
\De^2u &=0 \\
u|_{\dD} &=f_0\\
\f{\p u}{\p N}|_{\dD}&=\suml_{j=1}^{n} f_j N_j,\\
\norm{M(\nabla u)}{L^p(\dD)} &<\infty,
\end{array}\right.
$$
where $N_1,N_2, \ldots, N_{n}$ are the components of the normal vector and 
$f_0, f_1, f_2,\ldots, f_{n}$ satisfy the compatibility condition 
$(f_0,f_1,f_2, \ldots, f_n)\in WA_2(\dD)$ (cf. \cite{Pipher1}). 
The regularity problem is 
$$BR_p \left|\begin{array}{cc}
\De^2 u &=0 \\
D_2 u|_{\dD} &=f\\
\suml_{j=1}^{n-1}\dpr{\nabla_{T_j}}{\nabla D_j u}|_{\dD} &=g\\
\norm{M(\nabla\nabla u)}{L^p(\dD)} &<\infty
\end{array}\right.
$$
The $L^2$ theory for these problems (with the necessary 
adjustments for the order of the derivatives) is very similar to the 
harmonic case and we present it in Section \ref{sec:prel}. 
We have the following results in two dimensions.
\begin{theorem}
\label{theo:bdp}
There exists an $\ve=\ve(D)>0$, such that
if $0<\al<1/2+\ve$, $f_1,f_2\in C^\al(\dD)\cap L^2(\dD)$,  
then the unique $L^2$ solution to $BD_\al$ satisfies
$\nabla u\in C^\al(D)$. In fact,
$$
\norm{\nabla u}{C^{\al}(D)}+\supl_{X\in D}\textup{dist}(X,\dD)^{-1-\al}
|u(X)-u(X^*)-\dpr{X-X^*}{\nabla u(X^*)}|
\leq C\sum_{j=1}^2 \norm{f_j}{C^\al(D)},
$$
where $C$ is a constant depending only on the Lipschitz nature of $D$ and
$X^*$ is the projection of $X$ along the ``time'' axis  onto $\dD$.
Moreover the range $\al<1/2+\ve$ is sharp.
\end{theorem}
\begin{theorem}
\label{theo:brp} 
There exists an $\ve=\ve(D)>0$, such that the regularity problem $BR_p$ with
$2/3-\ve<p<2+\ve$, $(f,g)\in H^p_1(\dD)\times H^p(\dD)$ has unique (up 
to a constant) solution. Moreover the estimate
$$
\norm{M(\nabla\nabla u)}{L^p(\dD)}
\leq C \norm{\nabla_{T_1}f}{H^p(\dD)}+\norm{g}{H^p(\dD)},
$$
holds with a constant $C$ depending only on the Lipschitz nature of $D$.
The range $2/3-\ve<p$ is sharp.
\end{theorem}
\section{Preliminaries}
 \label{sec:prel}
 We separate this section into two parts - about harmonic and 
 biharmonic functions respectively. The corresponding equations 
 exhibit some common features like the $L^2$ theory, but there are some
 dissimilarities as well.  We try to present the similarities in the 
 technically simpler harmonic context and we briefly outline some 
 specifics for the biharmonic operator. We constantly refer in the 
 text to the papers \cite{Pipher2}, \cite{Pipher1} for the necessary 
 background results. Since the two dimensional case is of utmost 
 interest to us, we sometimes avoid the explicit formulas 
 (with the inevitable 
 technicalities that arise) for dimensions higher than two. 
 \subsection{Harmonic functions}
 \label{sec:Prel} Let $D\subseteq \rn$ be a Lipschitz domain, such
 that $D, D^{c}$ are connected. For technical reasons, we restrict
 our attention to the case of domains above Lipschitz graphs, i.e.
 \begin{eqnarray*}
  & & D=\{ (x,t): t>\varphi(x)\},\q \varphi:\mathbb{R}^{n-1}\to \rone,
\\
 & &|\varphi(x)-\varphi(y)|\leq M |x-y|.
 \end{eqnarray*}
 The surface measure on $\dD$ is defined via the usual
 $d\si=\sqrt{1+|\nabla \varphi|^2}dx$.

 Following \cite{Dahlberg1}, we introduce the atomic Hardy spaces
 $H^p(\dD)$ for $1\geq p>(n-1)/n$. First an $H^p(\dD)$ atom is a
 function $a:\dD\to \R$, such that
 \begin{eqnarray*}
 & & \textup{supp}(a)\subseteq B(Q,d)=\{P\in \dD: |P-Q|<d\}, \\ &
 &\int a d\si  =0, \\& &\norm{a}{L^2(\dD)}\leq C
 d^{(n-1)(1/2-1/p)}.
 \end{eqnarray*}
 Then,
 \begin{eqnarray*}& & H^p(\dD)=\{\sum \la_i a_i: \sum |\la_i|^p
 <\infty\}\\ & &\norm{f}{H^p(\dD)}=\ds\inf\limits_{f=\sum \la_i
 a_i} (\sum |\la_i|^p)^{1/p}.
  \end{eqnarray*}
 where $a_i$ are $H^p(\dD)$ atoms.

  We also define  $H^p_1(\dD)$
 atoms by requiring that \begin{eqnarray*} & &
 \textup{supp}(a)\subseteq B(Q,d)=\{P\in \dD: |P-Q|<d\}, \\ &
 &\norm{\nabla_T a}{L^2(\dD)}\leq C d^{(n-1)(1/2-1/p)},
 \end{eqnarray*}
 where $\ds\nabla_{T_j}u=\dpr{T_j}{\nabla u}=
 \left(\pp{}{x_j}+\pp{\varphi}{x_j}\pp{}{x_n}\right)u$.
  The  space $H^p_1(\dD)$ of
 distributions with one derivative in $H^p(\dD)$ may be defined as
 the $l^p$ span of such atoms. It is well known that
 $H^1(\dD)\subset L^1(\dD)$, while the spaces $H^p(\dD),\ p<1$
 contain non-integrable distributions.
Sometimes, we will abuse notations by writing $a(x)$, instead of 
$a(x,\vp(x))$. Observe that $\nabla_{T_j}a=\f{\p}{\p x_j}a(x,\vp(x))$.

We also define the (homogeneous) H\"older spaces $C^\al(\dD),\
 0<\al\leq 1$ by $$ C^\al(\dD)=\left\{f:\dD\to \rone:
 \norm{f}{C^\al(\dD)}=\sup\limits_{Q\neq P}
 \f{|f(P)-f(Q)|}{|Q-P|^\al}\right\}.$$ We remark
 that the H\"older spaces $C^\al$ and the Hardy spaces $H^p$, $\ds
 p=\f{n-1}{n-1+\al}$  can be paired in the sense that every element
 in one of them defines via integration a continuous linear
 functional on the other (cf. \cite{Stein} p. 130). Let
 $$\Gamma(x)=\left\{\begin{array}{cc}
 \ds\f{|x|^{2-n}}{(n-2)\omega_n} & n>2
 \\ \ds\f{1}{2\pi}\ln|x| & n=2\end{array}\right. $$ be the fundamental
 solution for the Laplace's equation in $\rn$. Define the single
 and double layer potentials  $\cs$ and $\ck$ by
\begin{eqnarray*}
 \cs(f)(X)  &=&  \pv \intl_{\dD} \Ga(X-Q) f(Q) d\si(Q),\q x\in
 \rn\setminus\dD \\
 \ck(f)(X)  &=&  \pv\intl_{\dD}\pp{\Gamma}{N_Q}(X-Q)f(Q)d\si (Q),\q
x\in
 \rn\setminus\dD
 \end{eqnarray*}
We also define the formal adjoint of $\ck$ $$ K^{*}(f)(X)=\pv
\intl_{\dD}\pp{\Gamma}{N_X}(X-Q)f(Q)d\si(Q).$$ 

\subsection{Biharmonic functions}
We start with the $L^2$ theory for the biharmonic equation, due to Kenig and
Verchota \cite{Kenigp} (see also Theorem 3.7 in \cite{Pipher2}).
\begin{proposition}
\label{prop:one}
Let $D\subset \rn$ be a Lipschitz domain. Then there exists an $\ve>0$ 
depending only on the Lipschitz nature of $D$, such that for 
$p:\q 2-\ve<p<2+\ve$
the equation 
$$\left|\begin{array}{cc}
\De^2u &=0 \\
u|_{\dD} &=f\\
\dpr{N}{\nabla u}|_{\dD} &=g\\
\norm{M(\nabla u)}{L^p(\dD)} &<\infty
\end{array}\right.
$$
is uniquely solvable. In addition, there are the estimates
\begin{itemize}
\item $|\nabla u(X)|\lesssim \textup{dist}(X,\dD)^{-(n-1)/p},$
\item $\norm{M(\nabla u)}{L^p\dD)}\lesssim \norm{\nabla u}{\dD}.$
\end{itemize}
\end{proposition}
There is a also the regularity result, which we now state( cf. Theorem 4.6,
\cite{Pipher2}).
\begin{proposition}
Let $D\subset \rn$ be a Lipschitz domain. Then there exists an $\ve>0$ 
depending only on the Lipschitz nature of $D$, such that for 
$p:\q 2-\ve<p<2+\ve$
the equation 
$$\left|\begin{array}{cc}
\De^2u &=0 \\
D_n u|_{\dD} &=f\\
\sum_{j=1}^{n-1}\dpr{\nabla_{T_j}}{\nabla D_j u}|_{\dD} &=g\\
\norm{M(\nabla\nabla u)}{L^p(\dD)} &<\infty
\end{array}\right.
$$
is uniquely solvable. In addition, there are the estimates
\begin{itemize}
\item $|\nabla\nabla u(X)|\lesssim \textup{dist}(X,\dD)^{-(n-1)/p},$
\item $\norm{M(\nabla\nabla u)}{L^p\dD)}\lesssim 
\sum_j(\norm{\nabla_{T_j}f}{L^p(\dD)}+\norm{g}{L^p(\dD)}).$
\end{itemize}
\end{proposition}
The fundamental solution of the biharmonic equation in two dimensions is
$$\Si(X)=\f{1}{8\pi} |X|^2\ln|X|.$$ Based on the $L^2$ theory, one is 
able to define the Green's function as follows.
Let $X$ be inside the domain and 
fix a point $X_0\notin\dD$. 
Let 
\begin{eqnarray*}
f_0(Q) &=& \Si(X-Q)-\Si(X_0-Q)=
\f{1}{8\pi}(|X-Q|^2 \ln|X-Q|-|X_0-Q|^2 \ln|X_0-Q|),\\
f_j(Q) &=& D_j(f_0(Q)),\q j\in \{1,2\}.
\end{eqnarray*}
Consider then the unique solution $\ga_X(Y)$ to the 
problem
$$
\left|\begin{array}{cc}
\De^2u &=0,\\
u|_{\dD} &=f_0,\\
\nabla u|_{\dD} &=\{f_1,f_2\}.
\end{array}\right.
$$
Since $\nabla_{T_j}f_j\sim |X-Q|^{-1}$, we have by the $L^2$ regularity
result that $\norm{M(\nabla\nabla \ga_X)}{L^2(\dD)}\lesssim 1$.
The Green's function can be defined as $G(X,Y)=\Si(X-Y)-\Si(X_0-Y)-\ga_X(Y)$.
Observe that since two tangential derivatives of the explicit quantity
$\Si(X-Y)-\Si(X_0-Y)$ also belong to $L^2(\dD)$, we can conclude
\begin{equation}
\label{eq:ltwo}
\norm{M(\nabla\nabla G(X,\cdot))}{L^2(\dD)}\lesssim 1.
\end{equation}
Later on, we will be able to show a much stronger estimate than 
\eqref{eq:ltwo} when the integration is over dyadic pieces away 
from the origin.

Integration by parts shows that one has the following representation formula 
for biharmonic functions (cf. \cite{Pipher1}):
\begin{equation}
\label{eq:ninetyone}
u(x)=\intl_{\dD}u(Q)\f{\p}{\p N_Q} \De_Q G(X,Q)+
\intl_{\dD}\f{\p u}{\p N_Q}(Q) \De_Q G(X,Q).
\end{equation}

\section {Existence and Uniqueness for harmonic functions}
\label{sec:existence}

 In this section, we prove existence and uniqueness for $R_p$ and $N_p$, and
Dirichlet problem with H\"older data in Theorem \ref{theo:main}. 
Our plan is as follows. First, we show
that the regularity and Neumann problems are equivalent, i.e. if
one can solve the regularity problem uniquely in $H^p_1(\dD)$,
then one can solve uniquely the Neumann problem in $H^p(\dD)$
and vice versa. Secondly, we show that the $R_p$ is solvable in
$H^p_1(\dD)$, $p>2/3-\ve$ and finally we prove uniqueness for
$R_p$. The uniqueness result will be almost automatic in view of
Lemma 2.2 in \cite{Brown} and the usual $L^q$ uniqueness result
for the Dirichlet problem, $q\geq 2-\ve$.
Finally, we show the existence and
uniqueness result for the Dirichlet problem with $C^\al(\p D)$
data in Theorem \ref{theo:th2}.
\subsection {Equivalence for the regularity and Neumann problem}

Let $u$ satisfy the Neumann problem $N_p$ with data $f\in \hpo$.
By the properties of $\hpo$, we consider without loss of
generality only smooth compactly supported data $f$. Let
$u_{-1}(x,t)=\intl_{0}^t u(x,z)dz$. Define $$v=a(x)
(u_{-1})_x+b(x) (u_{-1})_t.$$ It is not difficult to check that
 $$v_x(x,0)=(a (u_{-1})_x+b (u_{-1})_t)_x=-(b (u_{-1})_{xt}+
 c (u_{-1})_{tt})=-b u_x-c u_t=f(x).$$
We observe that $v$ satisfies $$ (R_p)\left| \begin{array}{cc}
\di(\tilde{A} \nabla v) &= 0 \\ v(x,0) &= g(x)
\end{array},\right.
$$ where $\tilde{A}=\f{A}{ac-b^2}$ and $g$ is the (unique)
function with $g'(x)=f(x), \int g=0$.

Since $A$ is an uniformly elliptic matrix with time independent
coefficients, so is $\tilde{A}$. Note that starting with a
solution of a Neumann problem, we have produced an $L^2$ solution
to an associated regularity problem. Moreover, since $\nabla
v=\left(\begin{array}{cc}-b & -c\\a & b
\end{array}\right)\nabla u$, one obtains pointwise equivalence $M(\nabla v)\sim
M(\nabla u)$. Hence, if one can prove estimates for $R_p$
$$\norm{M(\nabla v)}{L^p(\dD)}\lesssim \norm{g}{H^p_1(\dD)},$$
they would imply the corresponding estimates for $N_p$
$$\norm{M(\nabla u)}{L^p(\dD)}\lesssim \norm{f}{H^p(\dD)}.$$ Thus,
we have showed that if one can solve the regularity problem in
$H^p_1(\dD)$, then one can also solve the Neumann problem. The
reverse implication can be proved by retracing back the 
argument above, so we omit the details.

We note that the equivalence of the regularity and Neumann
problems in the sense described above is purely two dimensional
phenomena. Actually, in the important case of the Laplace's
equation, it is not difficult to check that $u$ and $v$ above are
in fact conjugate harmonic functions and thus one cannot expect the
equivalence to persist in higher dimensions.
Actually, by the equivalence of the regularity and Neumann problem 
and the existence results of Theorem \ref{theo:th2} (to be proved below), 
we establish Theorem \ref{prop:five}.

\subsection{Solvability for the regularity problem in $H^p_1(\dD)$}

By well known approximation techniques (see for example
\cite{Kenig}, section 1.10), it will suffice to prove the estimate
\begin{eqnarray}
\label{eq:regular} \norm{M(\nabla u)}{L^p}\leq C\norm{g}{H^p_1}
\end{eqnarray}
for solutions $u$ of  $R_p$ corresponding to smooth matrices $A$
and smooth data $g$, which are known to exist, as long as the
constant $C$ is independent of everything, but the ellipticity
constant of $A$.

Thanks to the ``atomic'' nature of $H^p_1(\rone)$, one can take
$g$ to be a $H^p_1$-atom. Indeed, if we assume \eqref{eq:regular}
for atoms and take into account the $p$-subaditivity of the
$L^p$ quasi-norm ($p<1$) we get for $g=\sum \la_i g_i$ 
$$
\norm{M(\nabla u_g)}{p}^p\lesssim \sum |\la_i|^p\norm{M(\nabla
u_{g_i})}{p}^p \lesssim \sum |\la_i|^p\lesssim \norm{g}{H^p_1}^p.
$$ Simple dilation and translation argument  allows us to reduce to the case
 of an unit atom, i.e.
\begin{enumerate}
\item $\textup{supp}\q g\subset [-1,1]$,
\item $\norm{g}{\infty}, \norm{g'}{\infty}\lesssim 1$.
\end{enumerate}
For $\tau\in (1/2,1)$, consider the intervals $R_j^\tau=(2^j
\tau,2^{j+1}/\tau)$ and $R_j:=R_j^1$. Let
$$q_j^\tau(x)=\left\{\begin{array}{ccc} & 100(2^j\tau -x)\qq &x<2^j
\tau,
\\ &0 \qq &2^j\tau\leq x\leq 2^{j+1}/\tau, \\ & 100(x-2^{j+1}/\tau), \qq
&x\geq 2^{j+1}/\tau,
\end{array}\right.
$$ and $\Omega_j^\tau=\left\{(x,t):t>q_j^\tau(x)\right\}$. Observe
that since $\Omega_j^\tau$ is a domain above Lipschitz graph, the
$L^2$ theory for divergence form equations with time independent
coefficients applies to it (see the discussion after Theorem
\ref{theo:main}). We have the following lemma.
\begin{lemma}
\label{lemma:klj}
 Let $u$ be the unique $L^2$ solution to the problem $$
\left|\begin{array}{cc} \di(A\nabla u) &=0\\ u(x,0) &=g(x),
\end{array}
\right. $$ where $g$ is an unit atom in $H^p_1(\rone)$. Then there exists
$\ve>0$, such that $$ \intl_{R_j}M(\nabla
u)^2\leq C_\ve 2^{(-\ve-2)j}.$$
\end{lemma}
Let us  show that Lemma \ref{lemma:klj} implies
\eqref{eq:regular}. By H\"older's inequality, $$
\intl_{R_j}M(\nabla u)^p\lesssim |R_j|^{1-p/2}
\left(\intl_{R_j}M(\nabla u)^2\right)^{p/2}\lesssim
2^{j(1-p/2+p(-\ve-2)/2)}, $$ and for every $p>2/3-O(\ve)$, 
the series $\sum_j
\int_{R_j}M(\nabla u)^p$ converges.

Thus, it remains to prove Lemma \ref{lemma:klj}.
\begin{proof}(Lemma \ref{lemma:klj})
We use the standard Cacciopoli type argument. By the $L^2$ theory
for $\Omega_j^\tau$ and since $R_j\subset
\partial\Omega_j^\tau$, we derive
\begin{equation}
\label{eq:123}
\intl_{R_j} M(\nabla u)^2\lesssim \intl_{1/2}^1\intl_{\p
\Omega_j^\tau}M(\nabla u)^2 d\tau\lesssim 2^{-j} \intl_{\Omega_j}
|\nabla u|^2,
\end{equation}
where $\Omega_j=\Omega_j^{1/2}$. Break $\Omega_j$
into ``good'' and ``bad'' part, so that
\begin{eqnarray*}
& &G_j=\Omega_j\cap\{t\geq 2^j\}\\ & & B_j=\Omega_j\setminus G_j.
\end{eqnarray*}
On the good part, we further decompose
$G_j=\bigcup\limits_{k=1}^{\infty} G_j^k$, so that
$G_j^k=G_j\cap\{t\sim 2^{j+k}\}$. For each $G_j^k$, one applies the
usual interior estimates for the solution (cf. \cite{Kenig}), to
get
\begin{equation}
\label{eq:good}
\intl_{G_j^k} |\nabla u|^2\leq 2^{-2(j+k)} \intl_{G_j^k}
|u|^2.
\end{equation}
By the $D_2$ solvability, we can always estimate 
$\norm{u}{L^2(G_j^k)}\lesssim 2^{(j+k)/2}$, which would give the 
desired estimate, except fot the extra decay factor $2^{-\ve j}$.\\
For the ``bad'' part,
select an even function  
$\psi\in \cc_0^{\infty}(\rtwo)$, 
so that supp $\psi\subset (1/4,4)\times(0,2)$ and  $\psi(x,t)=1$ 
for $(|x|,|t|)\in (1/2,2)\times(0,1)$. Denote 
$\psi_j:=\psi(2^{-j}\cdot)$.
Observe that since $u(x,0)=0$ on $R_j$, we may extend
$u(x,t)$ for $t<0$ across $R_j$ as an {\it even} function.
By the ellipticity of $A$ and 
the divergence theorem, one then derives
\begin{equation}
\label{eq:bad}
\intl_{B_j}|\nabla u|^2\leq \intl_{\rtwo}
\dpr{A\nabla (u\psi_j)}{\nabla (u \psi_j)}\lesssim 2^{-j} 
\intl_{\rtwo} |\nabla u| |u| \psi_j\lesssim 2^{-j}
\norm{\nabla u}{L^2(C_j)}\norm{u}{L^2(C_j)},
\end{equation}
where $C_j\supset B_j$ is again a box with 
sides $\sim 2^j$.

From \eqref{eq:bad}, $R_2$, $D_2$ solvability and by 
iterating \eqref{eq:bad} (we will get back to  this point later on),  
we easily get the 
bound $\int_{\Omega_j} |\nabla u|^2\leq C_{\de} 2^{-j}2^{\de j}$ for 
all positive $\de$. This estimate, together with  \eqref{eq:123} imply
Lemma \ref{lemma:klj} without the crucial term $2^{-\ve j}$.

The usual approach to get the improvement $2^{-\ve j}$  is to use
Sobolev embedding $H^1(\rn)\subset L^{2n/(n-2)}(\rn)$, which unfortunately 
{\it fails} for dimension two. 
We use instead the following multiplicative variant of Sobolev 
embedding
\begin{equation}
\label{eq:Sobol}
\norm{u}{L^4(\rtwo)}\leq 
\norm{u}{L^2(\rtwo)}^{1/2}\norm{\nabla u}{L^2(\rtwo)}^{1/2}.
\end{equation}
We have the following proposition.
\begin{proposition}
\label{cl:vajen} 
Suppose
\begin{eqnarray}
& &\norm{\nabla u}{L^2(Q_j)}^2\lesssim 
2^{-j} \norm{u}{L^2(P_j)}\norm{\nabla u}{L^2(P_j)},
\label{eq:1}\\
& &\norm{M(u)}{L^{2-\de}(\rone)}\lesssim 1,\label{eq:2}\\
& &\norm{M(\nabla u)}{L^2(\rone)}\lesssim 1 \label{eq:3}.
\end{eqnarray}
where $Q_j\subset P_j$ are boxes with sides $\sim 2^j$. Then 
 there exists $\ve=O(\de)>0$, such that 
\begin{equation}
\label{eq:crucial}
\norm{\nabla u}{L^2(Q_j)}\lesssim 2^{(-1/2-\ve)j}.
\end{equation}
\end{proposition}
It is clear that a direct application of \eqref{eq:crucial} gives 
the estimate for the ``bad'' part, while for the good 
part one applies \eqref{eq:crucial} for
$G_j^k$ and summation in $k>0$ then gives
\eqref{eq:good}. Hence, to complete the proof of 
Lemma \ref{lemma:klj}, it remains to prove Proposition \ref{cl:vajen}.
\end{proof}
\begin{proof}(Proposition \ref{cl:vajen})

Apply \eqref{eq:Sobol} for $u(x,t)\psi_j(x-x_0,t-t_0)$, where 
$(x_0,t_0)$ are suitably chosen so that 
$\psi_j(x-x_0,t-t_0)=1\q \textup{on}\q  P_j$ and supp 
$\psi_j(\cdot-x_0,\cdot-t_0)\subset 4 P_j$. Cauchy-Schwartz and 
\eqref{eq:2} yield 
\begin{eqnarray*}
& & \intl_{P_j} |u|^2\lesssim 
\norm{u}{L^{2-\de}(P_j)}\norm{u}{L^{(2-\de)'}(P_j)}
\lesssim 2^{j/(2-\de)-j/2)}2^{2j/(2-\de)'}\norm{u}{L^4(P_j)}\lesssim
\textup{by}\q \eqref{eq:Sobol} \\
& &\lesssim 2^{j(1/2+1/(2-\de)')}(\norm{u}{L^2(4 P_j)}^{1/2}
\norm{\nabla u}{L^2(4 P_j)}^{1/2}+
2^{-j/2}\norm{u}{L^2(4 P_j)})\lesssim \\
& & \lesssim 2^{j(1/2+1/(2-\de)')}(2^{j/4}
\norm{\nabla u}{L^2(4 P_j)}^{1/2}+1)\lesssim 
2^{(5/4-O(\de))j}\norm{\nabla u}{L^2(4 P_j)}^{1/2}+2^j2^{-O(\de)j},
\end{eqnarray*}
where $O(\de)$ is a positive number of the order of $\de$.
From the preceding estimate and \eqref{eq:1}, we get
\begin{equation}
\label{eq:wer}
\intl_{Q_j} |\nabla u|^2\lesssim 2^{-(3/8 +O(\de))j}
\norm{\nabla u}{L^2(4 P_j)}^{5/4}+2^{-j/2-O(\de)j}\norm{\nabla u}{L^2(4 P_j)}
\end{equation}
We can perform now the following iteration procedure. Call
$\la_j=\norm{\nabla u}{L^2(Q_j)}$ and $\mu_j=\norm{\nabla u}{L^2(4 P_j)}$.
Clearly \eqref{eq:wer} reads as
\begin{equation}
\label{eq:pol}
\la_j^2\lesssim 2^{-3/8 j} 2^{-O(\de)j}\mu_j^{5/4}+2^{-j/2}2^{-O(\de)j}\mu_j.
\end{equation}
Since \eqref{eq:3} allows us to bound $\mu_j\leq C 2^{j/2}$, one gets 
from \eqref{eq:pol} improvement for the bounds for $\la_j,\mu_j$. We continue 
in that fashion and use the improved bounds back at \eqref{eq:pol}. That way
one gets an improvement at every step.
One has
$$
\la_j,\mu_j \lesssim 2^{-j/2-\ve j} 
$$
for some $\ve=O(\de)$ and the proof is complete.
\end{proof}

\subsection{Uniqueness for the regularity problem in $H^1_p$}
The uniqueness result is almost automatic in two dimension due to 
the following lemma of Brown \cite{Brown}, which we state verbatim.
\begin{lemma}
\label{lemma:dve}
Let $D\subset \rn$ be a connected Lipschitz domain and suppose that $u$ is harmonic
in $D$. Let $X^*$ be a fixed point in $D$ and suppose $u(X^*)=0$. For 
$p<n-1$ and $p^*=(n-1)p/(n-1-p)$ we have
\begin{equation}
\label{eq:brov}
\norm{M(u)}{L^{p^*}(\p D)}\leq C\norm{M(\nabla u)}{L^p(\p D)},
\end{equation}
where the constant $C$ depends only on the distance of $X^*$ to the 
boundary, $p$ and the Lipschitz character of $\p D$.
\end{lemma}
Carefull inspection of the proof shows that one can relax the harmonicity
assumptions on $u$, by requiring that $u$ satisfies a divergence 
form equation. Indeed in the proof of \eqref{eq:brov}, Brown uses interior 
estimates and the equivalence of the square function with the non-tangential 
maximal function in $L^2$, which are of course available for solutions of
divergence form equations as well.

In contrast with the  higher dimensionional case, where
one needs to have an additional argument to 
prove uniqueness for $N_p$, $p>1-\ve$ (cf. \cite{Brown}), the two 
dimensional  uniqueness result follows from the
Brown's lemma for the range $1>p>2/3-\ve$ and  uniqueness
for $D_q$, $q>2-\ve$. To this end, 
assume that $u$ solves $R_p$ with zero data, such that 
$\norm{M(\nabla u)}{L^p(\rone)}<\infty$. From 
\eqref{eq:brov} we get
$$
\norm{M(u)}{L^{p^*}(\rone)}\leq C \norm{M(\nabla u)}{L^p(\rone)}.
$$
Observe that since $1>p>2/3-\ve$, we have $2-O(\ve)<p^*=p/(1-p)<\infty$ and therefore
$u$ solves a Dirichlet problem (with zero data), with 
$\norm{M(u)}{L^{p^*}(\rone)}<\infty$. Thus $u=0$ by the uniqueness for 
$D_{p^*}$.
\subsection{The Dirichlet problem with H\"older data}
We start with a lemma in the spirit of Theorem 3.4 in \cite{Brown}.
\begin{lemma}
\label{lemma:one} Let $D\subset \rtwo$ be a star-like Lipschitz
domain. There exist $\ve=\ve(D)>0$, so that for $2/3-\ve<p<1$ the
maps
\begin{eqnarray*}
\f{1}{2}I &+& K^{*}: \hpo\to \hpo \\ \f{1}{2}I &-& K^{*}: \hpo\to
\hpo
\end{eqnarray*}
are invertible.
\end{lemma}
Assuming the validity of Lemma \ref{lemma:one}, we can easily show
part three of Theorem \ref{theo:th2}.
Indeed, observe that $\f{1}{2}I +K :\ca\to\ca$ is the adjoint
map to $\f{1}{2}I +K^{*}: \hpo\to \hpo$ and $-\f{1}{2}I +K
:\ca\to\ca$ is the adjoint map to $-\f{1}{2}I +K^{*}: \hpo\to
\hpo$, where $\al=1/p-1$. Thus
\begin{eqnarray*}
\f{1}{2}I &+& K: \ca\to \ca \\ \f{1}{2}I &-& K: \ca\to \ca
\end{eqnarray*}
are invertible operators for $\al <
\left(\f{1}{2/3-\ve}-1\right)=1/2+O(\ve)$. Thus, the solution to
$$ (D_\al)\left| \begin{array}{cc} &\De u =0\\ & u|_{\dD} =f\in
C^\al(\dD)
\end{array}\right.$$ is in $C^{\al}(\bar{D})$.

{\bf Remark}
There exists more direct arguments
towards proving the $C^\al$ estimates of Theorem \ref{theo:th2}
similar to the one employed for the system of elastostatic
(\cite{Dahl2}), and later on for biharmonic functions
(\cite{Pipher1}).

The tools provided by Lemma \ref{lemma:one} however allow for
unified treatment of the problem at hand. More specifically, one
reduces the question for solvability of the regularity and Neumann
problems in  $H^p \q 2/3-\ve<p<2+\ve$ to the invertibility of
unitary perturbations of the adjoint of the double-layer potential
in the same spaces.

Lemma \ref{lemma:one}  follows from the existence and
uniqueness statements for $R_p$ and $N_p$, $2/3-\ve<p<1$ combined with
the usual duality argument.
The proof of Lemma \ref{lemma:one} is essentially contained in \cite{Brown} 
(cf. Proposition 3.1--3.5). One can easily adapt the argument there to the 
two dimensional case and the extended range of $p$'s, thus  we omit the 
details.

\section{Existence and Uniqueness for biharmonic functions}
\label{sec:exu}

In this section, we briefly sketch the proofs of Theorems 
\ref{theo:bdp}, \ref{theo:brp}. We follow closely the ad-hoc approach of
\cite{Pipher1}, which originated in \cite{Dahl2}. As we have mentioned
earlier, a more systematic way of studying the problem would be 
the  method  of Lemma \ref{lemma:one}, i.e. one could build an operator $T$,
whose invertibility is equivalent to the solvability of $BR_p$ and by duality
to the Dirichlet problem with $\al$-H\"older data in the sense of Theorem 
\ref{theo:bdp}. This program was implicitely carried out in \cite{Pipher2}.
We choose however the direct method for sake of clarity of the exposition.
\subsection{Uniqueness for biharmonic functions}
An easy adaptation of Lemma \ref{lemma:dve} 
gives the following.
\begin{lemma}
\label{lemma:three}
Suppose that for a given {\bf biharmonic} function $u$, 
there is $X^*\in D\subset \rn$, such that $|\nabla u(X^*)|=0$.
For $p<n-1$ and $p^*=(n-1)p/(n-1-p)$ there is
\begin{equation}
\label{eq:brown1}
\norm{M(\nabla u)}{L^{p^*}(\p D)}\leq C\norm{M(\nabla\nabla u)}{L^p(\p D)},
\end{equation}
where the constant $C$ depends only on the distance of $X^*$ to the 
boundary, $p$ and the Lipschitz character of $\p D$.
\end{lemma}
Indeed, one uses the equivalence of the area integral and 
the non-tangential maximal function for biharmonic functions 
as in the proof for the harmonic case to show \eqref{eq:brown1}. Since the 
maximum principle of \cite{Pipher1} is valid for dimensions two and three
(but not for dimensions bigger than three), we  argue as follows.
Assume that a biharmonic function solves $BR_p$ for $2/3-\ve<p<1$, such that 
$D_2 u|_{\dD}=0$ and  $\nabla_{T_1}D_1 u|_{\dD}=0$
and $\norm{M(\nabla\nabla u)}{L^p(\dD)}<\infty$.
Thus, after an eventual correction with a linear term, we may assume that
$\nabla u|_{\dD}=0$ and $u|_{\dD}=0$. By Lemma \ref{lemma:three}, we 
conclude
$$
\norm{M(\nabla u)}{L^{p^*}(\dD)}\lesssim 
\norm{M(\nabla\nabla u)}{L^p(\p D)}<\infty,
$$
where $2-\ve<p^*<\infty$. By the maximum principle, we have  uniqueness 
for the Dirichlet problem in $L^{p^*}$, hence $|\nabla u|=0$.
\subsection{Existence for the biharmonic regularuty problem}
For the existence part, we will use the following Cacciopolli type 
inequality.
\begin{lemma}
\label{lemma:greg}
Let $D\subset \rn$ be a domain above Lipschitz graph. Let 
$\Om_1\subset \Om_2\subset D$ be bounded Lipschitz domains. Let 
$\De^2u=0$ in $D$ with $M(\nabla^2u)\in L^2(\dD)$. 
Let $\ve$ is a small number as in Proposition \ref{prop:one}. Let
also $2-\ve<p<2+\ve$ and $d=\textup{dist}(\Om_1,D\setminus\Om_2)$.
Then there is a constant $C$, depending only on the Lipschitz constant and $p$,
so that 
\begin{eqnarray*}
\intl_{\Om_1}|\nabla^2 u|^2 dX &\leq& 
C(\norm{\nabla u}{L^{p'}(\dD\cap\p\Om_2)}\norm{M(\nabla^2 u)}{L^{p}(\dD)}+\\
&+& d^{-1} \norm{u}{L^{p'}(\dD\cap \p\Om_2)}\norm{M(\nabla^2 u)}{L^{p}(\dD)}+\\
&+& d^{-1} \norm{\nabla u}{L^2(\Om_2)}\norm{\nabla^2 u}{L^2(\Om_2)}+
d^{-2}\norm{u}{L^2(\Om_2)}\norm{\nabla^2 u}{L^2(\Om_2)}).
\end{eqnarray*}
\end{lemma}
Lemma \ref{lemma:greg} appears as Lemma 5.6 in \cite{Pipher2} for 
dimension three. The proof though can be easily adapted to this generality.
We state now our main estimate for ``atomic'' solutions.
\begin{lemma}
\label{lemmajk}
Let $a$ be an unit atom in $H^p(\dD)$, $D\subset \rtwo$, 
with $\int a(z)z=0$. There exists 
$\ve=\ve(D)>0$, such that the unique solution to the regularity problem 
$$BR_p \left|\begin{array}{cc}
\De^2 u &=0 \\
D_2 u|_{\dD} &=0\\
\dpr{\nabla_{T_1}}{\nabla D_1 u}|_{\dD} &=a\\
\norm{M(\nabla\nabla u)}{L^p(\dD)} &<\infty
\end{array}\right.
$$
satisfies 
\begin{equation}
\label{eq:four}
\intl_{\dD}M(\nabla^2u)^p\lesssim 1,
\end{equation}
for $2/3-\ve<p<1$.
\end{lemma}
{\bf Remark}
The extra cancellation condition $\int a(z)zdz=0$ is technical and it 
is  possible to remove. 
Simple translation and dilation argument yields 
 \eqref{eq:four} for arbitrary atoms $a$ in $H^p(\dD)$ with 
the special cancelation property $\int a(z)zdz=0$. 
Since such atoms  suffice to span $H^p(\dD)$, we get 
$$
\intl_{\dD}M(\nabla^2u_a)^p\lesssim \norm{a}{H^p(\dD)}^p
$$
In particular, we get \eqref{eq:four} for atoms without the 
extra cancellation $\int a(z)zdz=0$.

\begin{proof}(Lemma \ref{lemmajk})
We make the standard assumption that the boundary is smooth, so that 
smooth solutions
exist according to the classical theory. As usual, our estimates will 
not involve the smoothness constants and after one proves the result 
in that fashion, a standard approximation technique yields
\eqref{eq:four} for general Lipschitz domains.

Next, observe that the boundary conditions for $BR_p$ imply that 
$\f{d^2}{dx^2} u(x,\vp(x))=a(x)$ and therefore
$$
u(x,\vp(x))=\intl_{-\infty}^x\intl_{-\infty}^y a(z)dz dy.
$$
By support considerations and since $\int a=0$, 
$\int a(z)zdz=0$, we get
$u(x,\vp(x))=0$, for $x>2$ and
$u(x,\vp(x))=0$, for $x<-2$. Thus, $u$ also satisfies a {\it Dirichlet} type
boundary conditions
\begin{eqnarray*}
\left|\begin{array}{cc}
\De^2 u &=0 \\
D_2 u|_{\dD} &=0\\
u(x,\vp(x)) &=\intl_{-\infty}^x\intl_{-\infty}^y a(z)dz dy.
\end{array}\right.
\end{eqnarray*}
The advantage of casting $u$ as a solution to both Dirichlet and regularity
type problems will be seen  later on in the proof.

We first estimate \eqref{eq:four} for $x$-small. By H\"older and
$L^2$ regularity
$$
\left(\intl_{\dD\cap\{(x,\vp(x)):|x|<100\}}M(\nabla^2u)^p\right)^{1/p}
\lesssim \left(\intl_{\dD}M(\nabla^2 u)^2\right)^{1/2}\lesssim 1.
$$
For every point $X\in \dD$ fix a right cone $\Ga(X)$
opening upward with axis along the ``time'' axis and sides having slopes 
$100 \norm{\vp'}{\infty}$. Define the auxilliary maximal functions
\begin{eqnarray*}
M_1(\nabla^2 u)(X) &=\supl_{\Ga(X)\cap\Ga} |\nabla^2 u|,\\
M_2(\nabla^2 u)(X) &=\supl_{\Ga(X)\setminus\Ga} |\nabla^2 u|.
\end{eqnarray*}
$M_1$ measures the behavior of the solution away from the boundary and
is somewhat easier to control, while $M_2$ captures the behavior of
the solution close to the boundary. It is clear that $M\lesssim M_1+M_2$.

Since $D_2u=0$ and $\nabla u|_{\dD}\in L^{2-\ve}(\dD)$, we deduce from 
the $L^2$ Dirichlet theory 
$$
\norm{M(\nabla u)}{L^{2-\ve}(\dD)}\lesssim 1.
$$
Thus, interior estimates imply 
\begin{eqnarray*}
|\nabla u(X)| &\lesssim \textup{dist}(X,\dD)^{-1/(2-\ve)},\\
|\nabla\nabla u(X)| &\lesssim \textup{dist}(X,\dD)^{-1-1/(2-\ve)}.
\end{eqnarray*}
Therefore for $|Q|>100$, we infer
\begin{equation}
\label{eq:six}
M_1(\nabla^2 u)(Q)\lesssim |Q|^{-1-1/(2-\ve)}\in 
L^{2/3+O(\ve)}(\dD).
\end{equation}
For $R>10$ and $1\leq \tau\leq 2$ define the Carleson region $\Om_\tau^R$ 
above $Z_R=\{(x,\vp(x)):|x|\sim R\}$
$$
\Om_\tau=\Om_\tau^R=\left\{(x,t):R/\tau\leq |x|\leq R\tau,
\vp(x)<t<100\tau \norm{\vp'}{\infty}\right\}.
$$
By H\"older's inequality, we have
$$
\intl_{\dD\cap\p \Om_1 }M(\nabla^2u)^p
\lesssim R^{1-p/2}\left(\intl_{\dD\cap\p \Om_1}M(\nabla^2 u)^2\right)^{p/2}.
$$
Hence to show \eqref{eq:four} it will be enough to prove
$$
\intl_{\dD\cap\p \Om_1}M(\nabla^2 u)^2\lesssim R^{-2-\ve},
$$
for some positive $\ve>0$.

From the $L^2$ regularity result in $\Om_\tau^R$, we have
$$
\intl_{\p \Om_\tau\cap \dD}M_2(\nabla\nabla u)^2\lesssim  
\intl_{\p \Om_\tau\setminus\dD}|\nabla\nabla u|^2+
\intl_{\p \Om_\tau\cap \dD}|\nabla_{T_1}\nabla u|^2.
$$
Averaging in $\tau$ and  $\nabla_{T_1}\nabla u|_{\p\Om_t\cap\dD}=0$
yield
\begin{equation}
\label{eq:ten}
\intl_{\p \Om_1\cap \dD}M_2(\nabla\nabla u)^2\lesssim 
R^{-1}\intl_{\Om_2}|\nabla\nabla u|^2.
\end{equation}
The  boundary conditions $u|_{\p \Om_2\cap \dD}=0$, 
$\nabla u|_{\p \Om_2\cap \dD}=0$ and Lemma \ref{lemma:greg}  imply
\begin{equation}
\label{eq:seven}
\intl_{\Om_2}|\nabla^2 u|^2 \lesssim R^{-1} \norm{\nabla u}{L^2(\Om_3)}
\norm{\nabla^2 u}{L^2(\Om_3)}+
R^{-2}\norm{u}{L^2(\Om_3)}\norm{\nabla^2 u}{L^2(\Om_3)},
\end{equation}
for some eventually bigger box $\Om_3\subset D$ still 
having diameter $\sim R$. Since $u|_{\p \Om_3\cap \dD}=0$, one estimates
\begin{equation}
\label{eq:eleven}
\norm{u}{L^2(\Om_3)}\lesssim R\norm{\nabla u}{L^2(\Om_3)}.
\end{equation}
Combining \eqref{eq:ten},\eqref{eq:seven}, \eqref{eq:eleven} with the 
obvious $\norm{\nabla^2 u}{L^2(\Om_3)}\lesssim 
R^{1/2}\norm{M_2(\nabla^2 u)}{L^2(\p \Om_3\cap \dD)}$ yield
\begin{eqnarray}
& &\intl_{\Om_2}|\nabla^2 u|^2
\lesssim R^{-1/2} \norm{\nabla u}{L^2(\Om_3)}
\left(\intl_{\p \Om_3\cap \dD} M_2(\nabla^2 u)^2\right)^{1/2}
\label{eq:twenty} \\
\label{eq:eight}
& &\intl_{\p \Om_1\cap \dD}M_2(\nabla^2 u)^2 
\lesssim R^{-3/2} \norm{\nabla u}{L^2(\Om_3)}
\left(\intl_{\p \Om_3\cap \dD} M_2(\nabla^2 u)^2\right)^{1/2}
\end{eqnarray}
As usual, one has 
$\norm{\nabla u}{L^2(\Om_3)}\lesssim 
R^{1/2}\norm{M(\nabla u)}{L^2(\p \Om_3\cap \dD)}\lesssim R^{1/2}$ and by 
iterating \eqref{eq:eight}, one gets for every $\de>0$
\begin{equation}
\label{eq:forty}
\intl_{\p \Om_1\cap \dD}M_2(\nabla^2 u)^2 \leq C_\de R^{-2+\de},
\end{equation}
which barely fails to make $\intl_{\dD}M_2(\nabla^2 u)^p$ convergent. One 
also has
\begin{equation}
\label{eq:fifty}
\norm{\nabla^2 u}{L^2(\Om_1)}\lesssim 
R^{1/2}\norm{M_2(\nabla^2 u)}{L^2(\Om_1)}\leq C_\de R^{-1/2+\de/2}.
\end{equation}
Thus in order to get a better estimate we resort to \eqref{eq:Sobol}. 
Write 
\begin{eqnarray}
\label{eq:thirty}
\intl_{\Om_3} |\nabla u|^2 &\lesssim&
\left(\int_{\Om_3} |\nabla u|^{2-\ve}\right)^{1/(2-\ve)} 
\left(\int_{\Om_3} |\nabla u|^{(2-\ve)'}\right)^{1/(2-\ve)'} \\
\nonumber
&\lesssim& R^{1/2+1/(2-\ve)'} \norm{M(\nabla u)}{L^{2-\ve}(\dD)}
\left(\intl_{\Om_3}|\nabla u|^4\right)^{1/4}\q\textup{by}\q \eqref{eq:Sobol}\\
\nonumber
&\lesssim& R^{1-O(\ve)}\norm{\nabla u}{L^2(\Om_4)}^{1/2} 
\norm{\nabla\nabla u}{L^2(\Om_4)}^{1/2}+
R^{1/2-O(\ve)}\norm{\nabla u}{L^2(\Om_4)},
\end{eqnarray}
where $\Om_3\subset\Om_4$ is still a domain with diameter $\sim R$.
Since one can clearly derive \eqref{eq:forty}, \eqref{eq:fifty} for $\Om_4$
instead of $\Om_1$ (with eventually bigger constants), we use those with 
$\de=\ve/100$ to estimate the right hand side of \eqref{eq:thirty}. We get
\begin{equation}
\label{eq:sixty}
\intl_{\Om_3} |\nabla u|^2\lesssim 
R^{3/4-O(\ve)}\norm{\nabla u}{L^2(\Om_4)}^{1/2} +
R^{1/2-O(\ve)}\norm{\nabla u}{L^2(\Om_4)}.
\end{equation}
Iterate  \eqref{eq:sixty} to get 
$\norm{\nabla u}{L^2(\Om_3)}\lesssim R^{1/2-O(\ve)}$. Going back to 
\eqref{eq:eight} and using the newly obtained bound for 
$\norm{\nabla u}{L^2(\Om_3)}$, we have 
\begin{equation}
\label{eq:seventy}
\intl_{\p \Om_1\cap \dD}M_2(\nabla^2 u)^2 
\lesssim R^{-1-O(\ve)} 
\left(\intl_{\p \Om_3\cap \dD} M_2(\nabla^2 u)^2\right)^{1/2}.
\end{equation}
Iterating \eqref{eq:seventy} now gives the desired estimate
\begin{equation}
\label{eq:ninety}
\intl_{\p \Om_1\cap \dD}M_2(\nabla^2 u)^2\lesssim R^{-2-O(\ve)},
\end{equation}
which completes the proof of the Lemma.
\end{proof}
We finish the proof of Theorem \ref{theo:brp}, based on 
Lemma \ref{lemmajk}.
\begin{proof}(Theorem \ref{theo:brp})
Start with data $(f,g)\in H^p_1(\dD)\times H^p(\dD)$ in $BR_p$. Select a 
harmonic function $h$ with Dirichlet data $f$. Define $H=h_{-1}$ to be the 
primitive of $h$. Stein's lemma and the 
regularity estimates in Theorem \ref{theo:th2} imply
\begin{equation}
\label{eq:rop}
\norm{M(\nabla^2 H)}{L^p(\dD)}\lesssim \norm{M(\nabla D_2 H)}{L^p(\dD)}
\lesssim \norm{M(\nabla h)}{L^p(\dD)}\lesssim \norm{f}{H^p_1(\dD)}.
\end{equation}
By Theorem \ref{prop:five}
$$
\displaystyle\norm{\f{\p h}{\p N}}{H^p(\dD)}\lesssim \norm{f}{H^p_1(\dD)}.
$$
Lemma \ref{lemmajk} then provides a biharmonic function $v$, with data
$g-\f{\p h}{\p N}$ such that
\begin{equation}
\label{eq:ros}
\norm{M(\nabla^2 v)}{L^p(\dD)}\lesssim \norm{g-\f{\p h}{\p N}}{H^p(\dD)}. 
\end{equation}
Set $u=H+v$ to get a biharmonic function satisfying the desired boundary 
conditions. Combining \eqref{eq:rop}, \eqref{eq:ros} yields the estimate
$$
\norm{M(\nabla^2 u)}{L^p(\dD)}\lesssim \norm{(f,g)}{H^p_1(\dD)\times H^p(\dD)}.
$$
\end{proof}
\subsection{H\"older continuity of solutions}

We will show Theorem \ref{theo:bdp} based on the estimates for solutions 
of the regularity problem, in particular for the Green's function.
\begin{proof} (Theorem \ref{theo:bdp})
Recall the representation formula for biharmonic solutions 
\eqref{eq:ninetyone}
\begin{equation}
\label{eq:ninetytwo}
u(X)=\intl_{\dD}f_0(Q)\f{\p}{\p N_Q} \De_Q G(X,Q)d\si(Q)+
\intl_{\dD}\suml_{j=1}^2 f_j(Q) N_j(Q) \De_Q G(X,Q)d\si(Q).
\end{equation}
By rescaling and dilation, it will suffice to show
$$
|u(X)-u(X^*)-\dpr{X-X^*}{\nabla u(X^*)}|\lesssim \suml_{j=1}^2\norm{f_j}{C^\al(\dD)},
$$
when $dist(X,\dD)=1$.
To achieve that one obviously needs  estimates for  $\Delta G|_{\dD}$. 

Observe that by \eqref{eq:ninety}, \eqref{eq:crucial} and the construction of
the solution for the regularity problem (see the proof of Theorem 
\ref{theo:brp}), we can establish the following estimate for the 
Green's function of $\Delta^2$ in $D$
\begin{equation}
\label{eq:hundred}
\intl_{\dD\cap |Q-X^*|\sim 2^k} M(\nabla^2 G(X,Q))dQ\lesssim 2^{-2k}2^{-O(\ve)k},
\end{equation}
when $dist(X,\dD)=1$ (cf. (3.5) in \cite{Pipher1}).

For the first term in \eqref{eq:ninetytwo}, the normal derivatives of the
harmonic function $\Delta_Q G(X,\cdot)$ are converted into tangential 
derivatives for $f_0$ and one eventually replaces the harmonic function 
$\Delta_Q G(X,\cdot)$ by its
conjugate harmonic functions (cf. \cite{Pipher1}, p. 395). The estimates
that one needs are then the same ones as for the second term. We estimate
the second term in \eqref{eq:ninetytwo} below. 

Note that by considering 
$u(X)-u(X^*)-\dpr{X-X^*}{\nabla u(X^*)}$ instead of $u(X)$, we may assume 
without loss of generality that 
$f_0(X^*)=0$ and $f_j(X^*)=0$. Thus
\begin{eqnarray*}
& &\left|\intl_{\dD}\suml_j f_j(Q)N_j(Q)
\Delta_Q G(X,Q)d\si Q\right| = \\
&=&\left|\intl_{\dD}\suml_j(f_j(Q)-f_j(X^*))N_j(Q)
\Delta_Q G(X,Q)d\si Q\right|\\
&\lesssim& \suml_j \norm{f_j}{C^\al(\dD)}\intl_{\dD} |Q-X^*|^\al 
|\Delta_Q G(X,Q)|d\si(Q)\lesssim \\
&\lesssim& \suml_j \norm{f_j}{C^\al(\dD)}
\suml_{k=0}^{\infty}2^{k/2+\al k}\left(\intl_{|Q-X^*|\sim 2^k}|\Delta_Q G(X,Q)|^2
d\si(Q)\right)^{1/2}\lesssim\textup{by}\q\eqref{eq:hundred}\\
&\lesssim& \suml_j \norm{f_j}{C^\al(\dD)}\suml_{k>0} 2^{(\al-1/2-O(\ve))k}\lesssim
\suml_j \norm{f_j}{C^\al(\dD)},
\end{eqnarray*}
provided $\al<1/2+O(\ve)$.
\end{proof}
\section{Counterexamples}
\label{sec:Counter}
In this section, we provide counterexamples to show that the statement of
Theorem \ref{theo:th2} (and  consequently Theorem \ref{theo:main}) is sharp.
Since Theorems \ref{theo:bdp} and \ref{theo:brp} provide an extension and 
essentially contain the results of Theorem \ref{theo:th2}, 
the harmonic functions considered below should be viewed also as biharmonic 
counterexamples showing the sharpness  of the statements of 
Theorem \ref{theo:brp} and Theorem \ref{theo:bdp}.

Consider the domain 
$$
\Om_\de=\left\{z\in \cc: |z|<1,\q 0<\textup{arg}z<2\pi/(1+\de)
\right\},
$$
for some small $\de>0$. Observe that the domain has  Lipschitz constant
$O(1/\de)$ and is very ``non-convex'' as $\de\to 0$. 
We will show that the Dirichlet  problem with H\"older data in 
$C^\al(\p \Om_\de)$,  $\al>(1+\de)/2$ is not uniquely solvable.

Define the harmonic function 
$$u(z)=\textup{Im}\q z^{(1+\de)/2},
$$
with the obvious identification of the complex plane $\cc$ with $\rtwo$.
It is easy to check that $\p \Om_{\de}$ consist of two segments and an arc, so 
that $u$ vanishes on the segments and $u\in C^1$ on the arc. Altogether,
we have that $u$ solves a Dirichlet problem on $\Om_\de$ with  $C^1$ data, 
while one obviously cannot control more than $\norm{u}{C^{(1+\de)/2}(\Om_\de)}$.
Thus, we have shown that part three of Theorem \ref{theo:th2} is sharp.

Next, we invoke the equivalence results of Section 
\ref{sec:existence} to conclude that since one cannot solve the Dirichlet 
problem with H\"older data, then the regularity and Neumann problem must 
fail to be solvable as well.

\end{document}